\documentclass [12pt]{amsart}

\usepackage{amssymb,amsxtra,amsfonts}
\usepackage{graphics}

\usepackage{verbatim}  %for \comment
\usepackage{bbding}   % for \bigoasterisk----now use made-up version 

\newenvironment{ppb}[1]
{\ \!\!\!\!\!\!\!\!\!\!\!\!\!\!\!\!\!\!\!\!\!\!\!\!\!\!\!\!\!\!\!\!\!\!\!\!\!\!\!\! {\bf PPB------------------------------------------------------------------------------------------------PPB}\newline \tiny {#1}
\  \newline\normalsize\phantom{f}\!\!\!\!\!\!\!\!\!\!\!\!\!\!\!\!\!\!\!\!\!\!\!\!\!\!\!\!\!\!\!\!\!\!\!\!\!\!\!\! {\bf PPB------------------------------------------------------------------------------------------------PPB}\newline}{}

\def\reE@DeclareMathSymbol#1#2#3#4{%
    \let#1=\undefined
    \DeclareMathSymbol{#1}{#2}{#3}{#4}}
\DeclareSymbolFont{symbolsC}{U}{txsyc}{m}{n}
\SetSymbolFont{symbolsC}{bold}{U}{txsyc}{bx}{n}
\DeclareFontSubstitution{U}{txsyc}{m}{n}
\reE@DeclareMathSymbol{\strictiff}{\mathrel}{symbolsC}{76}

\newcommand\beq{\begin{equation}}
\newcommand\eeq{\end{equation}}
\newcommand\bal{\begin{align*}}
\newcommand\eal{\end{align*}}   %why does this not work??
\newcommand\bmx{\left(\begin{matrix}}
\newcommand\emx{\end{matrix}\right)}
\newcommand\bsmx{\left(\begin{smallmatrix}}
\newcommand\esmx{\end{smallmatrix}\right)}

\newcommand{\spq}{/\!\!/}

\newcommand{\st}{\ \bigl\vert\ }

\providecommand{\pdeg}{\text{\rm pdeg}}

\def\part#1{\frac{\partial\phantom{q}}{\partial#1}}

\newcommand {\flb}{\lbrack\!\lbrack}
\newcommand {\frb}{\rbrack\!\rbrack}
\newcommand {\flp}{(\!(}
\newcommand {\frp}{)\!)}

\newcommand{\MDR}{\mathcal{M}_{\text{\rm DR}}}
\newcommand{\MB}{\mathcal{M}_{\text{\rm B}}}

\newcommand{\MDol}{\mathcal{M}_{\text{\rm Dol}}}

\newcommand{\Jac}{\text{\rm Jac}}

\newcommand{\HH}{\text{\rm H}}

\newcommand{\Lie}{{\mathop{\rm Lie}}}
\newcommand{\rank}{\mathop{\rm rank}}

\newcommand{\res}{{\mathop{\rm Res}}}

\newcommand{\tr}{{\mathop{\rm Tr}}}

\DeclareMathOperator{\Hom}{Hom}         % this looks to be the correct way to do this (put a * ``...otar*{'' if want things underneath

\newcommand{\GL}{{\mathop{\rm GL}}}

\newcommand{\End}{\mathop{\rm End}}

\newcommand{\hk}{{hyperk\"ahler }}   %useful word
\newcommand{\IB}{\mathbb{B}}
\newcommand{\IC}{\mathbb{C}}

\newcommand{\IH}{\mathbb{H}}

\newcommand{\IP}{\mathbb{P}}                                     
                           
\newcommand{\IR}{\mathbb{R}}

\newcommand{\IZ}{\mathbb{Z}}

\newcommand{\cD}{\mathcal{D}}
\newcommand{\cE}{\mathcal{E}}

\newcommand{\cF}{\mathcal{F}}

\newcommand{\cK}{\mathcal{K}}

\newcommand{\M}{\mathcal{M}}
\newcommand{\cM}{\mathcal{M}}

\newcommand{\cO}{\mathcal{O}}

\newcommand{\g}{       \mathfrak{g}     }
\newcommand{\gk}{      \mathfrak{g}_k   }

\newcommand{\gks}{     \mathfrak{g}_k^*   }
\newcommand{\gkis}{    \mathfrak{g}_{k_i}^* }

\newcommand{\lt}{\mathfrak{t}}
\newcommand{\lh}{\mathfrak{h}}

\newcommand{\lp}{\mathfrak{p}}

\renewcommand{\ll}{\mathfrak{l}}

\newcommand{\wh}{\widehat}

\newcommand{\al}{\alpha}

\newcommand{\be}{\beta}

\newcommand{\Ga}{\Gamma}

\newcommand{\La}{\Lambda}

\newcommand{\si}{\sigma}

\newcommand{\Si}{\Sigma}
\renewcommand{\th}{\theta}

\makeatletter
 \newlength{\typesize}
\makeatother

\newlength{\vvoff}
\newlength{\hhoff}

\newcommand{\pf}{\begin{bpf}}

\newcommand{\pfms}{\begin{bpfms}}
\newcommand{\epf}{\end{bpf}\hfill$\square$\\}           % end proof
\newcommand{\epfms}{\end{bpfms}\hfill$\square$\\}       % end proof

\newcommand{\idea}{\begin{bidea}}

\newcommand{\eidea}{\end{bidea}\hfill$\square$\\}           % end proof

\newcommand{\sk}{\begin{bsk}}    %type: \sk ..... \esk

\newcommand{\esk}{\end{bsk}\hfill$\square$\\}           % end sketch
\newcommand{\sketch}{\begin{bsketch}}%type: \sketch ..... \esketch

\newcommand{\esketch}{\end{bsketch}\hfill$\square$\\}

\newtheorem {hypo}{\bf\hspace{-\parindent}Hypothesis}
\newtheorem {thm}[hypo]{Theorem}   %remove [hypo] to number separately
\newtheorem {prop}[hypo]{Proposition}%[section]
\newtheorem {lem}[hypo]{Lemma}%[section]

\newtheorem {defn}[hypo]{Definition}%[section]
\theoremstyle{remark}\newtheorem{rmk}[hypo]{Remark}

\newenvironment{exercise}{{\bf\hspace{-\parindent}Exercise.}}{}

\begin{document}

\title[Nonabelian {H}odge theory of irregular curves]{Hyperk\"ahler manifolds and nonabelian {H}odge theory of (irregular) curves}
\author{Philip Boalch}%

\begin{abstract}
Short survey based on talk given at the Institut Henri Poincar\'e January 17th 2012, during program on surface groups. The aim was to describe some background results before describing in detail (in subsequent talks) the results of \cite{gbs} %
related to wild character varieties and irregular mapping class groups.
\end{abstract}

\maketitle

\section{Big picture}

Lets start by recalling the usual picture for nonabelian Hodge theory on curves, due to Hitchin, Donaldson, Corlette and Simpson 
\cite{Hit-sde, Don87, Cor88, Sim-hbls}.

Fix an integer $n$ and let $G=\GL_n(\IC)$.
Let $\Si$ be a smooth compact complex algebraic curve.
Given this data one may consider the nonabelian cohomology space
$$\cM = \HH^1(\Si,G).$$
Ignoring stability conditions for the moment (until the next section), this space is naturally a {\em hyperk\"ahler  manifold}, and there
are three viewpoints on it:

1) (Dolbeault)  as the moduli space $\MDol$ of Higgs bundles, consisting of pairs  $(E, \Phi)$ with $E\to \Si$ a rank $n$ degree zero holomorphic vector bundle and $\Phi\in \Ga(\End(E)\otimes \Omega^1)$ a Higgs field,

2) (De Rham) as a moduli space $\MDR$ of connections on rank $n$ holomorphic vector bundles, consisting of pairs  $(V, \nabla)$  with  $\nabla:V\to V\otimes \Omega^1$ a holomorphic connection, and

3) (Betti) as the space $\MB=\Hom(\pi_1(\Si),G)/G$ of conjugacy classes of representation of the fundamental group of $\Si$.

This gives three different algebraic structures on the same underlying space $\cM$  (since $\Si$ is compact, by GAGA, the holomorphic objects above are in fact algebraic, and have algebraic moduli spaces).
$\MDR$ and $\MB$ are complex analytically isomorphic via the Riemann--Hilbert correspondence, taking a connection to its monodromy representation.
$\MDR$ and $\MDol$ are naturally diffeomorphic as real manifolds via the nonabelian Hodge correspondence,   but are not complex analytically isomorphic. 
Thus there is more than one natural complex structure on $\cM$; they form part of the family of complex structures making $\cM$ into a \hk manifold.

To get an idea of this first consider the abelian case $n=1$, so $G=\IC^*$.

Then one finds:

$\bullet$ $\MDol\cong T^*\Jac(\Si)$ is the cotangent bundle
of the Jacobian variety of $\Si$,

$\bullet$ $\MDR\to \Jac(\Si)$ is a twisted cotangent bundle
of the Jacobian of $\Si$; it is an affine bundle modelled on the cotangent bundle,

$\bullet$ $\MB\cong (\IC^*)^{2g}$ is isomorphic to $2g$ copies of $\IC^*$.

In this case one may compute explicitly the Riemann--Hilbert isomorphism 
$\MDR \to \MB$
and see it involves exponentials and so is not algebraic. A slightly stronger statement is also true:

\begin{lem}
\footnote{Hint: use the valuative criterion for properness for the compositions 
$\IC^* \hookrightarrow \MB\to \MDR\to \Jac(\Si)$.}
There is no algebraic isomorphism $\MB \to \MDR$.
\end{lem} 

On the other hand since $\MB$ is affine it has no compact holomorphic 
subvarieties, and so neither does $\MDR$.
Thus it is clear there is no complex analytic isomorphism 
$\MDR\to \MDol$ since the zero section $\Jac(\Si)\to \MDol$ is a compact holomorphic subvariety.
On the other hand using (abelian) Hodge theory and the Dolbeault isomorphisms it is easy to obtain a (non-holomorphic) isomorphism $\MDR\to \MDol$.

In this abelian case the corresponding hyperk\"ahler metric is {\em flat}, but in general it is highly nontrivial, and difficult to make explicit.

Returning to the general picture, the three different viewpoints have different applications, and are of interest to different groups of people for different reasons. For example:

1) The Dolbeault spaces are algebraically completely integrable Hamiltonain systems (the Hitchin systems): there is a proper map 
$$\MDol\to \IH$$
to a vector space of half the dimension, the generic fibres of which are abelian varieties. (In the abelian case the space was a product $\IC^g\times \Jac(\Si)$, but in general the fibres vary nontrivially and there are singular fibres),
 
2) The Betti spaces are the complex character varieties of $\Si$, and as such they admit a natural (symplectic, algebraic) action of the mapping class group of $\Si$, coming from the natural action of the mapping class group on the fundamental group $\pi_1(\Si)$, 

3) When the curve $\Si$ varies in a family over a base $\IB$ the corresponding De Rham spaces assemble in to a fibre bundle over $\IB$, which has a natural flat algebraic (Ehresmann) connection on it: the nonabelian Gauss--Manin connection. 
When written in explicit coordinates this gives a natural class of nonlinear differential equations coming from geometry, generalising both the classical isomonodromy/Painlev\'e VI equations, and the usual (abelian) Gauss--Manin connections, abstracting the (linear) Picard--Fuchs equations. Integrating this nonlinear connection around a loop in $\IB$ gives a transcendental automorphism of the fibre $\MDR$: upon conjugating by the Riemann--Hilbert map these give  {\em algebraic} automorphisms of the Betti spaces as in 2).
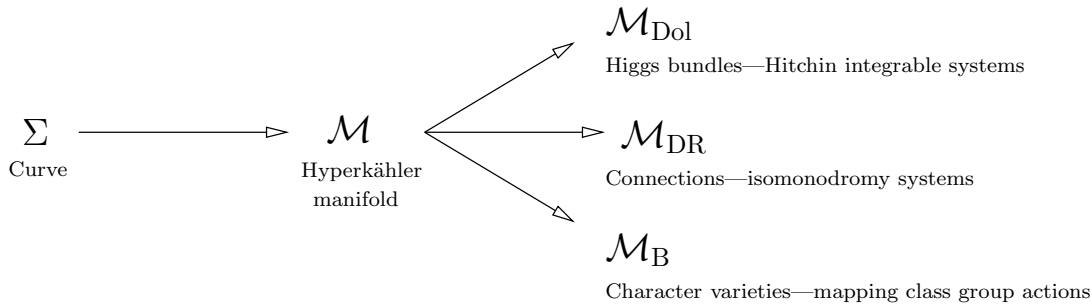
\begin{figure}[h]
\input{scheme.pstex_t}
\end{figure}

This is a very rich picture which many people like for many (usually different) reasons.
One natural extension is to replace the initial curve by a higher dimensional projective variety: this has been done (mainly by Simpson) and has many important applications. 
But the moduli spaces that arise are smaller than the spaces that occur for a curve: in brief if one chooses a sufficiently  generic curve in the variety the restriction map (restricting say a flat connection to the curve) yields an injective map embedding the moduli space into that for the curve (cf. Fujiki \cite{fujiki-hk} p.3).

Remaining with curves then, the next generalisation is to consider punctured curves, and to consider parabolic structures at the punctures. In essence one considers connections/Higgs fields  with simple poles and compatible with the parabolic structures.
In this context the nonabelian Hodge correspondence was established by Simpson \cite{Sim-hboncc} and such \hk metrics were constructed  by Konno \cite{Kon} and Nakajima \cite{Nak}. 
In this context the algebraic integrable systems and the isomonodromy equations have a long history in the case when $\Si$ is the Riemann sphere.

But one can obtain many more moduli spaces by considering meromorphic connections with higher order poles: the ``wisdom'' (if one can call it that) garnered by studying this case (and the corresponding isomonodromy equations and the Fourier--Laplace transform) is that one should consider the extra parameters that control the  coefficients of the connections/Higgs fields  which are more singular than the residue, as being {\em analogous to the moduli of the curve}.

Indeed one obtains {\em new} discrete group actions extending  the usual braid/mapping class group actions in 2) above by varying these extra parameters. This motivates the following definition.
\subsection{Irregular curves}

The basic aim is to replace the initial curve $\Si$ in the above story by an ``irregular curve'' defined as follows.
Some more motivation and examples will appear in the next section.
Fix the group $G$ and a maximal torus $T\subset G$, and denote the Lie algebras
$\lt\subset \g$.

\begin{defn}\label{defn: irreg curve}
An ``irregular curve'' consists of 

1) a smooth compact complex algebraic curve  $\Si$, 

2) distinct marked points $a_1,\ldots,a_m\in \Si$, and

3) an irregular type $Q_i$ at $a_i$ for $i=1,\ldots,m$.
\end{defn}

In turn an `irregular type' is defined as follows.
(We will give the coordinate independent definition first, and then explain it in coordinates---the abstract viewpoint will be useful later when we vary the curve.)
Let $\wh \cO_i$ denote the formal completion of the ring of germs at $a_i$ of holomorphic functions on $\Si$, and let $\wh \cK_i$ denote its field of fractions.

\begin{defn} \label{def: irtype}
An ``irregular type'' $Q_i$ at $a_i$ is an element
$$Q_i \in \lt(\wh \cK_i) /\lt(\wh \cO_i).$$ 
\end{defn}

One may think of an irregular type as a $\lt$-valued meromorphic function germ, well defined modulo holomorphic terms. 
Explicitly, if we choose a  local coordinate $z$ on $\Si$ vanishing at $a_i$,
then 
$\wh\cO_i= \IC\flb z \frb, \wh \cK_i = \IC\flp z \frp$, and 
 $Q_i$ may be  written in the form
$$Q_i = \frac{A_{r_i}}{z^{r_i}} + \cdots + \frac{A_{1}}{z}$$
for elements $A_j\in \lt$ for $j=1,\ldots, r_i$, for some $r_i> 0$.

Given an irregular type $Q_i$ as above there is a well defined subgroup $H_i\subset G$ (the centralizer of $Q_i$)
defined as $H_i=\{g\in G\st gA_jg^{-1}=A_j, j=1,\ldots,r_i\}$;
it is a reductive subgroup of $G$ again with maximal torus $T$, and we will
write $\lh_i=\Lie(H_i)$. 

Now we wish to attach moduli spaces $\MDol$ and $\MDR$ to an irregular curve (the Betti picture will be the focus of the subsequent three lectures).
As in the punctured case we also need to use parabolic structures, so will first discuss flags/filtrations.

\subsection{Flags/filtrations}\label{ss: flags}

Let $\lt_\IR=X_*(T)\otimes\IR\subset \lt$ be the real cocharacters
(the real span of the lattice $X_*(T)$ of one parameter subgroups).

Given $\th\in \lt_\IR$ there is a canonically determined parabolic subgroup
\beq\label{eq: parab}
P_\th = \{ g\in G\st z^\th g z^{-\th} 
\text{ has a limit as $z\to0$ along any ray}\}
\eeq
with Lie algebra $\lp_\th=\Lie(P_\th)$.
Since we are working here with $\GL_n(\IC)$ it is easy to see what this means in terms of matrix entries. Moreover we can reinterpret it in terms of real filtrations: For any $\al\in\IR$ let $E_\al\subset V:=\IC^n$ be the 
$\al$-eigenspace of $\th\in \End(V)$, and define a filtration $\cF_\th$ of $V$ as follows:
$$(\cF_\th)_\be = \bigoplus_{\al\ge \be} E_\al\subset V$$
for any $\be\in \IR$.
Then %
$P_\th$ is just the subgroup of $G$ preserving the filtration $\cF_\th$.
(Note it will be important for us not to insist that we are working in a basis for which the diagonal entries of $\th$ are ordered.)

\subsection{Irregular connections}

Fix an irregular curve $\Si = (\Si, \{a_i\},\{Q_i\})$ and a weight 
$\th_i\in \lt_\IR$ at each marked point.
We will assume that each diagonal entry of each weight satisfies
\beq\label{eq: wt condn}
0\le (\th_i)_{jj} <1.\eeq
Here, for $G=\GL_n(\IC)$, this is no loss of generality and simplifies the presentation (avoiding discussion of parahoric structures/filtered $\cD$-modules), 
but for other reductive groups one does not get the full picture without considering parahoric structures (see \cite{logahoric}).

Let $k_i=r_i+1$, where $r_i$ is the order of the pole of $Q_i$, and define the divisor $D=\sum k_i(a_i) $ on $\Si$.
We will consider the moduli space of triples $(V,\nabla,\cF)$
where

$\bullet$ $V\to \Si$ is a rank $n$ holomorphic vector bundle,

$\bullet$ $\nabla:V\to V\otimes \Omega^1(D)$ is a meromorphic connection on $V$ with poles bounded by $D$,

$\bullet$ $\cF$ consists of a filtration $\cF_i$ of the fibre $V_{a_i}$ of $V$ at $a_i$ for each $i=1,\ldots, m$.

These data should be such that near each point $a_i$ there is a local trivialization of $V$ such that 

1) $\nabla = d-A$ where 
$$A = dQ_i + \La_i \frac{dz}{z} + holomorphic\  terms$$
for some $\La_i\in \lh_i$ in the Lie algebra of the centralizer of $Q_i$, where $z$ is a local coordinate vanishing at $a_i$, 

2) the filtration $\cF_i$ equals 
the standard filtration $\cF_{\th_i}$ on $\IC^n$ determined by the weight $\th_i$ (using the isomorphism $V_{a_i}\cong \IC^n$ coming from the trivialization),

3) the residue $\La_i$ preserves the filtration $\cF_i$ (i.e. 
$\La_i\in \lp_{\th_i}$).

Finally in order to obtain complex symplectic rather than Poisson moduli spaces (which is essential if we want to construct \hk manifolds!) we need to fix a certain adjoint orbit, as follows.

The conditions imply $\La_i$ is in the Lie subalgebra 
$\lh_i\cap \lp_{\th_i}\subset \lh_i$. 
This is just the parabolic subalgebra 
of $\lh_i$ determined by $\th_i$ (defined as in \eqref{eq: parab}, replacing 
$G$ by $H_i$).
Thus there is a projection
$$\pi:\lh_i\cap \lp_{\th_i} \to \ll_i$$
from this parabolic onto its Levi factor $\ll_i$ (quotienting by the nilradical).
Since we have fixed $\th_i$ we may identify $\ll_i$ with the centralizer of $\th_i$ in $\lh_i$:
$$\ll_i = \{ X\in \lh_i\st [\th_i,X] = 0\}.$$
It is again a reductive Lie algebra with Cartan subalgebra $\lt$.
The extra data to be fixed is the adjoint 
orbit of 
$$\pi(\La_i)\in \ll_i.$$
Using the Jordan decomposition any such orbit contains an element of the form
$$\tau_i+\si_i +N_i$$
where $\tau_i+\si_i \in \lt$ with ($\tau_i\in \lt_\IR, \si_i\in \sqrt{-1}\lt_\IR$) and $N_i\in \ll_i$ a nilpotent element commuting with $\tau_i+\si_i$.

Given such a meromorphic connection $V=(V,\nabla,\cF)$ its parabolic degree is:
$$\pdeg(V) = \deg(V) + \sum_1^m\tr(\th_i)$$
where $\deg(V)\in \IZ$ denotes the  degree of the underlying vector bundle $V$ on $\Si$. 

The notion of subconnection is defined in the usual way, as a subbundle $U\subset V$ preserved by $\nabla$, with the induced connection and filtrations in the fibres over the marked points.

A connection $V$ is {\em stable} if for any subconnection $U$ one has
$$\frac{\pdeg(U)}{\rank{U}} < \frac{\pdeg(V)}{\rank{V}}.$$
It is {\em semistable} if the weaker condition ($\le$) holds.

Let $\MDR(\Si, \th,\tau,\si,N)$
denote the moduli space of isomorphism classes of stable meromorphic connections $(V,\nabla,\cF)$ on the irregular curve $\Si$ 
with parabolic degree zero, weights $\th=\{\th_i\}$ and orbits 
determined by $\tau=\{\tau_i\}, \si=\{\si_i\}$ and $N=\{N_i\}$.

\begin{rmk}
It is known that any meromorphic connection on a disc is meromorphically isomorphic to one of the above form (with irregular part $dQ$ diagonal), after possibly passing to a finite cover (taking a root of the local coordinate $z$), so this restriction is not as strong as it might seem. Moreover the methods used here extend directly to the general case involving such ramifications.
\end{rmk}

On a curve a connection with a simple pole is the same thing as a logarithmic connection; we will henceforth refer to the case when all $Q_i=0$, as the logarithmic case. (Recall that a ``regular singular connection'' on $\Si\setminus\{a_i\}$  is a connection on an algebraic vector bundle on $\Si\setminus\{a_i\}$, that is isomorphic to the restriction of a logarithmic connection on a vector bundle on  
$\Si$.)

\subsection{Irregular Higgs bundles}

Similarly one may consider Higgs bundles on irregular curves, as follows.

Fix an irregular curve $\Si = (\Si, \{a_i\},\{Q_i\})$ as before and choose a weight $\th'_i\in \lt_\IR$ at each marked point
with $0\le (\th'_i)_{jj} <1$.

We will consider the moduli space of triples $(E,\Phi,\cE)$
where

$\bullet$ $E\to \Si$ is a rank $n$ holomorphic vector bundle,

$\bullet$ $\Phi:E\to E\otimes \Omega^1(D)$ is a meromorphic Higgs field with poles bounded by $D$,

$\bullet$ $\cE$ consists of a filtration $\cE_i$ of the fibre $E_{a_i}$ of $E$ at $a_i$ for each $i=1,\ldots, m$.

These data should be such that near each point $a_i$ there is a local trivialization of $E$ such that 

1) in this trivialisation the Higgs field takes the form
$$\Phi = -\frac{dQ_i}{2} + \Ga_i \frac{dz}{z} + holomorphic\  terms$$
for some $\Ga_i\in \lh_i$, 

2) the filtration $\cE_i$ equals 
the standard filtration $\cE_i=\cF_{\th'_i}$ on $\IC^n$ determined by the weight $\th'_i$,

3) the residue $\Ga_i$ preserves the filtration $\cE_i$.

Again we fix the adjoint orbit of $\Ga_i$ under the 
projection
$\pi:\lh_i\cap \lp_{\th'_i} \to \ll'_i$
where
$\ll'_i = \{ X\in \lh_i\st [\th'_i,X] = 0\}$
is the centralizer in $\lh_i$ of $\th'_i$.
We parameterise such orbits by choosing 
$\tau'_i+\si'_i +N'_i\in \ll_i'$
where $\tau'_i+\si'_i \in \lt$ has real part $\tau_i'$ and 
$N'_i\in \ll'_i$ is nilpotent and commutes with $\tau'_i+\si'_i$.

Given such a meromorphic Higgs bundle $E=(E,\Phi,\cE)$ its parabolic degree is:
$$\pdeg(E) = \deg(E) + \sum_1^m\tr(\th'_i)$$
and the notion of sub-Higgs bundle is defined in the usual way, as a subbundle $U\subset V$ preserved by $\Phi$, with the induced Higgs field and filtrations at the marked points.
A Higgs bundle $E$ is {\em stable} if for any sub-Higgs bundle $U$ one has
${\pdeg(U)}/{\rank{U}} < {\pdeg(E)}/{\rank{E}}.$

Let $\MDol(\Si, \th',\tau',\si',N')$
denote the moduli space of isomorphism classes of stable meromorphic Higgs bundles $(E,\Phi,\cE)$ on the irregular curve $\Si$ 
with parabolic degree zero, weights $\th'$ and orbits 
determined by $\tau', \si'$ and $N'$.

The main result of \cite{wnabh} %
can then be stated as 
\begin{thm}\label{thm: main}
The moduli space $\MDR(\Si,\th,\tau,\si, N)$ 
of meromorphic connections is a \hk manifold and is
 naturally diffeomorphic to the moduli space
\newline
\noindent
$\MDol(\Si,\th',\tau',\si', N')$
of meromorphic Higgs bundles if $N=N'$ and the parameters are related by:
$$\th_i' = -\tau_i-[-\tau_i],\qquad \tau_i' = -(\tau_i+\th_i)/2,\qquad
\si_i' = -\si_i/2$$
where $[\,\cdot\,]$ denotes the (component-wise) integer part.
Moreover the \hk metrics are {\em complete} if $N=0$ and there are no strictly semistable objects, and this may be ensured by taking the parameters to be off of some explicit hyperplanes.
\end{thm}

If all the irregular types $Q_i$ are zero this correspondence 
reduces to the tame case 
established by Simpson \cite{Sim-hboncc}, and
the `rotation' of the parameters 
$\th,\tau,\si$ is essentially the same.
If there are no marked points $m=0$ this reduces to the  original case of compact Riemann surfaces (of  Hitchin, Donaldson, Corlette and Simpson).
The general correspondence is again established by passing through solutions to Hitchin's self-duality equations, and the map from meromorphic connections to such solutions 
(i.e. the irregular analogue of the result of Donaldson \cite{Don87} and Corlette \cite{Cor88} constructing a harmonic metric for irregular connections on curves) was established earlier by Sabbah \cite{Sab99}.
As usual Hitchin's equations constitute the vanishing of a \hk moment map, and the \hk quotient of \cite{wnabh} is a strengthening of the infinite dimensional complex symplectic quotient description of $\MDR$ in the irregular case from \cite{Boa, smid}, generalising the Atiyah--Bott approach \cite{AB83}.

As mentioned in \cite{wnabh} many spaces of meromorphic Higgs bundles have been shown to be algebraically completely integrable Hamiltonian systems (generalised Hitchin systems) by Bottacin and Markman in \cite{ Bot, Mar},
extending  Hitchin \cite{Hit-sbis} in the holomorphic case  and e.g.  Adams et al \cite{AHH} and Beauville \cite{ Bea} in genus zero.
(The properness of the Hitchin map for such meromorphic Higgs bundles was established by Nitsure \cite{Nit-higgs}.)
Thus Theorem \ref{thm: main} extends the class of algebraic integrable systems  known to admit natural \hk metrics on their total space (although we seem to be working in slightly more generality due to the parabolic structures at the irregular singularities).

\begin{exercise}
Show $N_i\in \g$ commutes with $\th_i,\tau_i$ if and only if it commutes with $\th_i',\tau'_i$.
\end{exercise}

\section{More motivation}

The irregular connections will have fundamental solutions which, near a marked point $a_i$ with irregular type $Q_i$, involve essentially singular terms of the form $\exp(Q_i)$. 
We will give some basic examples of situations were such connections arise.

\subsection{Fourier--Laplace}
One often encounters the idea that irregular connections should be viewed as ``degenerations'' of logarithmic connections, as singular points coalesce. Whilst this viewpoint may be useful for some purposes, 
it is hard to relate the moduli spaces which appear before and after the coalescence.

Rather, the first example of an irregular connection to have in mind should probably  be that which arises by taking the Fourier--Laplace transform of a logarithmic connection on the Riemann sphere.
The basic fact that the function  $\exp(z)$ satisfies the differential equation
$$df = f dz$$
which has a pole of order two at $z=\infty$, implies that the Fourier--Laplace transform of a solution of a logarithmic connection will be a solution of an irregular connection (see \cite{BJL81} for a basic class of examples). And moreover by considering stability conditions etc this can be shown to induce an algebraic {isomorphism} between the corresponding De Rham moduli spaces (and between the corresponding Betti spaces). 
In such cases the irregular connections that arise have only two poles on the Riemann sphere, a pole of order one and a pole of order two. 
Thus the simplest irregular moduli spaces are in fact isomorphic to the logarithmic case (so the moduli spaces are certainly not more ``degenerate'' in the irregular case---this is also reflected in the fact that irregular moduli spaces admit {\em complete} hyperk\"ahler metrics.) In general of course there will be many moduli spaces (even on the Riemann sphere), that are not isomorphic to a logarithmic case.

Another aspect of this example of the Fourier--Laplace transform, that motivated the definition of irregular curve, is the fact that on the logarithmic side the moduli of the curve basically amounts to the positions of the poles (and we know that their motion underlies many of the known braid group actions), whereas on the isomorphic irregular side there are only two poles (so have no moduli), but there is a nontrivial irregular type at the pole of order two: the choice of the pole positions on the logarithmic side matches up precisely with the choice of the irregular type on the irregular side, i.e. the moduli of the (irregular) curves on each side matches up.
(Moreover the isomonodromy equations/braiding etc match up---this may be understood in terms of Harnad's duality \cite{Harn94} which can be shown to be equivalent to Fourier--Laplace (cf. \cite{k2p}\S3, \cite{yamakawa-mc+hd, slims}).

\subsection{Baker functions}
The Krichever construction \cite{Krich-methods} gives explicit algebro-geometric solutions to many soliton equations, and the key ingredient is the 
Baker-Akhiezer function: a function on an auxhiliary algebraic curve (the spectral curve) which is meromorphic except at some marked points where it has a prescribed essential singularity of the form $\exp(x/z+t_2/z^2+t_3/z^3+\cdots)$ in terms of a local coordinate $z$ vanishing at the marked point (beware in soliton theory one often takes the local coordinate to have a pole at the marked point, i.e. replaces $z$ by $z^{-1}$): such functions are horizontal sections of an irregular connection on a line bundle on the spectral curve (and the ``times'' $t_1=x,t_2,t_3,\ldots$ correspond to the choice of irregular type). See also for example \cite{schilling-bakerfns, SW}.
Thus in the abelian case $G=\IC^*$ the deformations of the irregular curve correspond to the ``times'' appearing in integrable hierarchies.

\subsection{Nonabelian Hodge theoretic invariants}
A standard way to find invariants of algebraic varieties is via Hodge theory. For example the map taking the ray through the cohomology class of the holomorphic symplectic form on a marked K3 surface gives a locally injective map to a period domain  classifying such surfaces.
More generally there are examples where one attaches to a variety a Hodge structure of an associated variety (cf. e.g.  \cite{ACT-cubic.3fds}).

The notion of Frobenius manifold was introduced \cite{Dub95long} to axiomatize 2d topological quantum field theories (TQFTs), extending Atiyah's axioms \cite{Ati89} (which for dimension two say that a TQFT is a Frobenius algebra) to include the natural deformations as well 
(so a Frobenius manifold is a family of Frobenius algebras with certain extra properties, encoding solutions of the WDVV equations of Witten et al).
In any case Dubrovin \cite{Dub95long} gave a local classification of the class of {\em semisimple} Frobenius manifolds: in essence one attaches a irregular differential system on the Riemann sphere, now called the {\em quantum differential equation}, to the Frobenius manifold, and takes the isomorphism class of that. In other words we take a point of a nonabelian cohomology space of a variety 
($\IP^1$) attached to the Frobenius manifold. And this gives a locally injective map on the moduli space. (Globally one needs to quotient by the braid group action coming from isomonodromy.) Thus, whereas usual Hodge theory may be used to classify classical objects such as algebraic varieites, nonabelian Hodge theory may be used to classify more complicated objects such as TQFTs.

\section{Examples of moduli spaces}

\subsection{Additive approximations}
Whilst in the usual case of holomorphic connections the moduli spaces attached to the Riemann sphere are trivial, in the meromorphic case there are many highly nontrivial moduli spaces with $\Si=\IP^1$.
On the Riemann sphere one can get a good idea of the moduli spaces by considering the open subset
$$\cM^*\subset \MDR$$
consisting of connections on bundles which are globally holomorphically trivial.
(Ignoring stability  this is an open subset of a component of $\MDR$.)
One can describe $\cM^*$ explicitly  as a finite dimensional complex symplectic quotient as follows (here we assume the weights $\th_i=0$):

Recall  $k_i= r_i+1$, where $r_i$ is the pole order of $Q_i,$ for 
$ i=1,\ldots,m$
and the connections have poles $D=\sum_1^m k_i(a_i)$.
Consider the group 
$G_{k}$ of $k$-jets of bundle automorphisms:
$$G_{k}:= \GL_n\bigl(\IC[z]/z^{k}\bigr).$$
The Lie algebra $\gk$ of $G_k$ consists of elements
$$X=X_0+X_1z+\cdots+X_{k-1}z^{k-1}$$
with $X_i\in \g=\Lie(G)$, and the Lie bracket is given in the obvious way, truncating terms of order $z^k$ or more.
One can identify the dual $\gks$ with the set of elements
\begin{equation} 	\label{eqn: pp map}
A=A_{k}\frac{dz}{z^{k}}+\cdots+ A_1\frac{dz}{z}
\end{equation}
with $A_i\in \g$, where we pair such  element $A$ with $X\in \g_k$ via the pairing
$$\langle A,X\rangle := \res_0(\tr(A\cdot X))=
\sum_{i=1}^{k}\tr(A_iX_{i-1}).$$

Thus, being the dual of a Lie algebra, $\gks$ is naturally a Poisson manifold
and its symplectic leaves are the coadjoint orbits.
On the other hand it is natural to identify a point of $\gks$ with the polar part of a meromorphic connection having a pole of order $k$, on the trivial bundle over a disc with local coordinate $z$.  
For each $i=1,\ldots,m$ let
$$\cO_i\subset \g_{k_i}^*$$
be the coadjoint orbit through the point
$$dQ_i + \La_i\frac{dz}{z}\in \g_{k_i}^*$$
taking $z$ to be a local coordinate vanishing at $a_i$, for 
$\La_i=\tau_i+\si_i+N_i\in \lh_i\subset \g$.

\begin{prop}(cf. \cite{smid} \S2)
Ignoring stability conditions, the space $\M^*(\Si, 0,\tau,\si,N)$ is isomorphic to the complex symplectic quotient
\begin{equation} \label{eqn: mod sp}
 \cO_1\times\cdots\times \cO_m\spq G\quad
\end{equation}
of the product of these coadjoint orbits $\cO_i\subset \gkis$ by the diagonal action of the constant group $G\subset G_{k_i}$ at the zero value of the moment map.
\end{prop}

For example if all the irregular types are zero, so $G_{k_i}=G$ for all $i$, and $\cO_i\subset \g^*\cong \g$, then $\cM^*$ is just 
 $$\cO_1\times\cdots\times \cO_m\spq G = 
\{(A_1,\ldots,A_m)\st A_i\in \cO_i, \sum A_i=0 \}/G$$
and a point of $\cM^*$ corresponds to the isomorphism class of a Fuchsian system of the form
$$d-\sum\frac{A_i}{z-a_i}dz$$
so that $\cM^*$ is just the moduli space of Fuchsian systems (with fixed residue orbits).
Specialising further, if $G=\GL_2(\IC)$ these spaces are much studied
under the name ``hyperpolygon spaces''.

Note that many of the spaces $\cM^*$ also have complete hyperk\"ahler metrics (for example in the Fuchsian case since Kronheimer, Biquard and Kovalev showed that the complex coadjoint orbits have invariant hyperk\"ahler metrics, and one may interpret the complex symplectic quotient above as a \hk quotient by the unitary subgroup of $G$---in fact here, for $G=\GL_n(\IC)$, the Fuchsian moduli spaces also arise as finite dimensional hyperk\"ahler quotients via Nakajima's work on quiver varieties). 
But these metrics on $\cM^*$ are {\em different} to the restriction of the full metric on $\cM$: in general both metrics are complete,  
but they are on different spaces ($\cM$  is a partial compactification of $\cM^*$: there really are stable connections on nontrivial degree zero bundles, cf. \cite{wnabh} Lemma 8.3).

\subsection{Low dimensional examples}
One of the advantages of studying the general meromorphic case is that there are nontrivial examples of moduli spaces of complex dimension two, i.e. real dimension four: thus the moduli spaces are hyperk\"ahler four-manifolds, and as such are ``gravitational instantons'' in the physics terminology. 
From a purely mathematical viewpoint Atiyah \cite{atiyah-hk}\S3 has emphasised that hyperk\"ahler four manifolds deserve special attention since they are the quaternionic analogues of Riemann surfaces or algebraic curves.

Most examples of the spaces $\cM$ of dimension two were first detected in the form of nonlinear differential equations: when written in explicit coordinates the corresponding irregular nonabelian Gauss--Manin connection amounts to a {\em second order} nonlinear differential equation. The classical Painlev\'e equations  arise in this way, and so in effect one may now translate the list of Painlev\'e equations into a list of complete hyperk\"ahler manifolds, and in turn there are the corresponding algebraically integrable (Hitchin) systems---which in this dimension are certain rational elliptic fibrations.

The list of known examples of dimension two is basically the following:
$$
\boxed{
\begin{matrix}
E_8 & E_7 & E_6 &      &          &  D_2 & (D_1) & (D_0) \\
    &     &     &  D_4 &  A_3 &      &      &    \\   
    &     &     &      &          & A_2  & A_1  & (A_0)   \\   
\end{matrix}
}
$$
$$\text{Table 1}$$

The symbols represent the affine Weyl symmetry groups that arise (mainly through work of Okamoto \cite{OkaP24, OkaPVI, OkaPV, Okamoto-dynkin} and Sakai \cite{Sakai-CMP01} in the context of Painlev\'e equations); the subscript (multiplied by three) is the real dimension of the space of parameters $\th,\tau,\si$ that arises.
All of the spaces to the right of and including $D_4$ arise in the theory of Painlev\'e differential equations, as follows:
\begin{center}
  \begin{tabular}{| c | c | c |c | c | c | c | c| c| }
\hline
{Space} 
&   $D_4$ &  $A_3$ &  $A_2$  & $A_1$  & $A_0$ & $D_2$ & $D_1$ & $D_0$ \\
\hline
{Painlev\'e equation} 
&   $6$ &  $5$ &  $4$  & $2$  & $1$ & $3$ & $3'$ & $3''$ \\
\hline
  \end{tabular}
\end{center}
where $3'$ and $3''$ denote special cases of the third Painlev\'e equation \cite{OKSO}.

A key point to understand is that each of these spaces arises as a moduli space of connections in many different ways: use of the Fourier--Laplace transform for connections on the Riemann sphere yields many isomorphisms between moduli spaces. Nonetheless in the literature there are certain standard representations (where the rank of the vector bundles is minimal)---this is usually expressed in terms of finding a ``Lax pair'' for the nonlinear differential equations.
For the cases corresponding to the Painlev\'e equations 
the standard representations of these spaces are as follows: one takes $G=\GL_2(\IC)$ and considers connections on rank two bundles over the Riemann sphere with the following pole orders:
\begin{center}
 \begin{tabular}{| c | c | c |c | c | c | c | c| c| }
\hline
{Space} 
&   $D_4$ &  $A_3$ &  $A_2$  & $A_1$  & $(A_0)$ & $D_2$ & $(D_1)$ & $(D_0)$ \\
\hline
{Pole orders} 
&   $1111$ &  $211$ &  $31$  & $4$  & $4$ & $22$ & $22$ & $22$ \\
\hline
  \end{tabular}
\end{center}

(The spaces in parentheses have some nilpotent leading terms and so require passing to a two-fold cover to diagonalise the irregular type---the argument of \cite{wnabh} shows they too are complete; in fact the nilpotents terms imply there can be no nontrivial subconnections in these three cases.)
Thus for example the Painlev\'e VI case corresponds to connections on rank two bundles over $\IP^1$ with four first order poles, and the others may be viewed as ``coalescences'' of this case.
(In the context of parabolic Higgs bundles this case was pointed out as the simplest nontrivial case much more recently by Boden--Yokogawa \cite{BodYok}.)
Most of these realisations are well-known (at least to the Painlev\'e community), and appear as Lax pairs in \cite{JM81, OKSO} (classically they were written in terms of second order differential operators, rather than connections on rank two vector bundles).
In particular for the spaces to the right of (and including) $A_3$ the only known realisations are as moduli spaces of irregular connections whereas the other four spaces ($D_4$ and the $E$s) do arise as spaces of logarithmic connections.

The spaces $E_6,E_7,E_8$ only admit higher rank realisations however, for example as logarithmic connections having three poles on the Riemann sphere on bundles of ranks $3,4,$ or $6$ respectively, and certain choices of residue orbits. They are related to Painlev\'e {\em difference} equations  cf. \cite{quad} and references therein, but not to nonlinear {\em differential} equations; the underlying curve has no nontrivial deformations: in the other cases the {\em irregular} curve has a one dimensional space of nontrivial deformations, yielding the `time' in the Painlev\'e equation. One can read off the Lax pairs from the corresponding affine Dynkin graphs: for example for $E_8$

\begin{figure}[h]
\input{dynkinE8.pstex_t}

\text{Affine $E_8$ Dynkin diagram}
\end{figure}
\noindent
one may take logarithmic connections with three poles on rank $6$ bundles, and generically the three orbits at the poles have respectively:
1) 6 distinct eigenvalues, 
2) 3 eigenvalues each repeated twice, and 
3) 2 eigenvalues each repeated three times (corresponding to the three legs of the affine $E_8$ graph).

It seems natural to conjecture that, if we add the cotangent bundle of an  elliptic curve (with flat metric), then up to isomorphism all deformation classes of two-dimensional \hk wild nonabelian Hodge moduli spaces appear in Table 1 (noting that there are ``twists'' that occur by varying the topological type of the underlying vector bundle).

My present understanding is that the spaces $D_i$, $i=0,1,2,3$, first appeared as hyperk\"ahler manifolds in the work of Cherkis--Kapustin \cite{chka1, chka2} (this includes  $A_3=D_3$), 
whereas the logarithmic examples are essentially covered by the general construction of Konno and Nakajima \cite{Kon, Nak}, 
leaving the cases of $A_2,A_1$, corresponding to the second and fourth Painlev\'e equations, as new four-dimensional examples covered by the general construction of \cite{wnabh}\footnote{
The fact that Painlev\'e moduli spaces are gravitational instantons has perhaps not been well absorbed by the community, although the author has been stating this in talks for years---most conspicuously perhaps \cite{iastalk07}, after which Witten enquired about the $A_0$ case.}.
(A different approach to such \hk four-manifolds
has been given recently by Hein \cite{hein-jams}.)

In all these cases, apart from Painlev\'e 3, the open subset 
$\cM^*\subset \cM$ is isomorphic to the asymptotically locally Euclidean (ALE) hyperk\"ahler manifold denoted by the same symbol
(cf. \cite{quad} exercise 3):
The ALE spaces exist for any simply laced affine Dynkin graph (i.e. of type $A_n,D_n,E_6,E_7$ or $E_8$)
and were constructed in general by Kronheimer \cite{Kron.ale} 
(those of type $A$ are due to Gibbons--Hawking \cite{gib-haw78}, which specialises to the  Eguchi--Hanson space 
for $A_1$). See e.g. Atiyah \cite{atiyah-hk}\S3 or Hitchin \cite{hitchin-bbk} for more discussion of these ALE spaces; the underlying spaces are deformations of the minimal resolution of the corresponding Kleinian surface singularity.
(For $\cM$ of type $D_2$, i.e. Painlev\'e 3, the space $\cM^*$ is isomorphic to the $D_2$ asymptotically locally flat (ALF) space.)
Thus we see a picture where the spaces $\cM$ appear as ``more transcendental versions'' of some quite well-known hyperk\"ahler manifolds: the ALE spaces arise as finite dimensional \hk quotients whereas we only know how to construct the \hk manifolds $\cM$  as \hk quotients of infinite dimensional spaces.

\subsection{$E_8$ examples}
Finally  we will describe the spaces corresponding to the case of $E_8$ in more detail (for generic parameters).
See \cite{quad} and references therein for more discussion of points 1-3) below and \cite{eor-gda} for 4).

Let $C\subset \IP^2$ be a cuspidal cubic curve, and let $\IC\subset C$ be its smooth locus (isomorphic to the affine line, hence the notation).
Choose $9$ points $p_1,\ldots,p_9\in \IC$, and blow-up $\IP^2$ at these $9$ points in order, and then remove the strict transform of $C$, leaving a  
noncompact complex surface $S$. 

1) If the $9$ points add up to zero (in the group law on $\IC$, which can be taken to coincide with that from the cubic---zero should be the unique point of inflection) then $S\cong \MDol(E_8)$ is isomorphic to one of the Dolbeault spaces of type $E_8$. 
The $9$ points are then basepoints of a pencil of cubics, and so $S$ is fibred by elliptic curves (the Hitchin fibration in this context).
In particular if all $p_i=0$ then one is repeatedly blowing up the inflection point and the tangent line decouples on the third blow up, leaving an affine $E_8$ arrangement of $(-2)$ curves in the interior of $S$.

2) If the $9$ points do not add up to zero, e.g. if the points in 1)  are perturbed a little bit, then $S\cong \MDR(E_8)$ is isomorphic to one of the De Rham moduli spaces (by scaling, the nonzero value of the sum is irrelevant and can be set to be $1\in\IC$).

3) If one repeats the story but with only $8$ points, one obtains one of the corresponding additive spaces $\cM^*(E_8)$ (in this case if all the points are zero one obtains an $E_8$ arrangement of $(-2)$ curves in the interior, as one would expect on the resolution of an $E_8$ Kleinian singularity),

4) Finally the Betti space looks completely different and one would probably not have guessed they were analytically isomorphic to the De Rham spaces without being told: here one repeats the above story but with eight points on the smooth locus $\cong \IC^*$ of a {\em nodal} cubic in $\IP^2$ (see \cite{eor-gda}).

\section{Conclusion}

Thus we have shown that it is possible to replace the curve in the diagram on p.3 with an ``irregular curve'' and still obtain associated \hk manifolds, which we have described from the Dolbeault and De Rham points of view, as moduli spaces of meromorphic Higgs bundles and connections respectively. 
It seems clear that most deformation classes of complete \hk manifolds of any given dimension arising from Higgs bundles on curves, arise only in the irregular case (as we saw in Table 1 in the case of real dimension 4).
To obtain a more explicit description of the underlying differentiable manifolds we will describe the Betti approach in the following three talks, and its extension to arbitrary complex reductive groups $G$  (following \cite{gbs}). In this description, for generic parameters, the moduli spaces arise as smooth affine algebraic varieties, defined by explicit equations. Further we will describe an algebraic way to obtain their holomorphic symplectic structures, and discuss the notion of admissible deformations of an irregular curve, which leads to a generalisation of the well-known mapping class group actions on the character varieties.

\renewcommand{\baselinestretch}{1}              %
\normalsize
\bibliographystyle{alpha}    \label{biby}
\bibliography{../thesis/syr} 

{\small
\'Ecole Normale Sup\'erieure et CNRS, 
45 rue d'Ulm, 
75005 Paris, 
 France

www.math.ens.fr/$\sim$boalch\ \ \ \ \ \ \qquad 

boalch@dma.ens.fr 

\nopagebreak
A version of this text appears at: http:/\!/hal.archives-ouvertes.fr/hal-00672481}

\end{document}